\documentclass[oneside,draft,12pt]{amsart}
\usepackage{amssymb, amscd}
\usepackage[all]{xy}

\theoremstyle{plain}
\newtheorem{theorem}{Theorem}[section]

\newtheorem*{theorem*}{}
\newtheorem{proposition}[theorem]{Proposition}
\newtheorem{corollary}[theorem]{Corollary}

\theoremstyle{definition}
\newtheorem{conjecture}[theorem]{Conjecture}
\newtheorem{definition}[theorem]{Definition}
\newtheorem{question}[theorem]{Question}
\newtheorem{example}[theorem]{Example}
\newtheorem{remark}[theorem]{Remark}

\DeclareMathOperator{\Tr}{tr}

\begin{document}

\title{Constructing and Deconstructing Group Actions}
\vskip .5in
\author{Alejandro Adem }

\address{Mathematics Department\\
         University of Wisconsin\\
         Madison, Wisconsin, 53706} 
\email{adem@math.wisc.edu}
\thanks{Partially supported by the NSF. This paper is a modified version 
of the plenary
lecture delivered by the author at the March 2002 Algebraic Topology Conference 
held at Northwestern University}
\vskip .5in

\begin{abstract}
Given a finite group $G$, it is not hard to show that it can
act freely on a product of spheres. A more delicate issue
is the following
question:
what is the minimum integer $k=k(G)$ such that $G$ acts
freely on a finite complex with the homotopy type of
a product of $k$ spheres?
The study of this problem breaks up into two distinct aspects:
(1)~proving bounds on $k(G)$ in terms of subgroup data
for $G$; and 
(2)~constructing explicit actions on a product of $k$ spheres.
In this note we will discuss aspects of both problems, including
recent progress based on the periodicity methods developed in \cite{ASm}.
We also
describe a 
potential counterexample to the prevalent expectations as well
as recent work on constructing geometric actions on actual
products of spheres.

\end{abstract}

\maketitle
\vskip .5in

\section{Introduction}

In 1944 P. A. Smith proved \cite{Sm}
that if a finite group $G$ acts freely on a sphere, then all of its abelian
subgroups must be cyclic.
This condition is known to be equivalent to the \emph{periodicity} of
of the cohomology of $G$ with trivial integral coefficients--recall that
the cohomology of $G$ is said to be periodic if there exists 
a $d>0$ such that 
$H^i(G,\mathbb Z)\cong H^{i+d}(G,\mathbb Z)$
for all $i>0$.
These groups are characterized by having generalized quaternion
$2$--Sylow subgroups and cyclic $p$--Sylow subgroups for $p$ odd.

In 1957, Milnor \cite{Mi}
showed that there are restrictions on $G$ imposed by
the geometry of the action:
if $G$ acts freely on a sphere $\mathbb S^n$, then every involution
in $G$ must be central. For example, the symmetric group on three
letters $\Sigma_3$ cannot act freely on any sphere, even though it has
periodic cohomology.
In 1960 Swan \cite{S} showed that a converse of Smith's result
does hold for homotopy spheres; he proved that 
if $G$ has periodic cohomology then it acts freely on a finite
complex $X\simeq \mathbb S^n$ for some $n>0$.
After that there ensued substantial activity on the problem of spherical
space forms, which culminated in 1976 with the solution 
due to
Madsen, Thomas and Wall \cite{MTW}; namely they proved that
a finite group $G$ acts freely on some sphere $\mathbb S^n$
if and only if all of its abelian subgroups are cyclic and
every element of order 2 in $G$ is central.
In other words, the combination of Smith's homological
condition and Milnor's geometric condition were enough to
ensure the existence of a free action on a sphere.

Based on the above one can inquire about the analogous situation
for groups which do not have periodic cohomology. In particular
given such a group, can we construct a free action on a 
suitably \emph{minimal} product
of spheres? 
In this note we will survey aspects of this problem, and in particular
describe recent substantial progress for groups which do not
contain $\mathbb Z/p\times\mathbb Z/p\times\mathbb Z/p$ as
a subgroup.

\section{Free Actions on a Product of Spheres}

We can define
the $p$--rank of a finite group $G$ as
$r_p(G) =~max~\{ r~|~ (\mathbb Z/p)^r\subseteq G \}$
and its rank  
as $r(G)=~max~\{r_p(G)~|~p~\rm{divides}~ |G|\}.$
Smith's pioneering work can be reformulated as saying
that if a finite
group $G$ acts freely on a finite complex 
$X\simeq \mathbb S^n$, then $r(G)=1$. Motivated by this and
other evidence,
the following rather difficult conjecture emerges:

\begin{conjecture}\label{conj-A}
If $G$ acts freely on
$X=\mathbb S^{n_1}\times\dots\times \mathbb S^{n_m}$, then
$r(G)\le m$.
\end{conjecture}

This question has been settled in the affirmative 
for $m=2$ by
Conner (in the equidimensional case \cite{Co}) and 
by Heller \cite{H}.
Better general 
results exist for an arbitrary product of equidimensional
spheres; we recall a result due to Adem and Browder
\cite{AB}.

\begin{theorem}
If $G$ acts freely on a finite complex
$X\simeq (\mathbb S^n)^k$, then
$r_p(G)\le k$ for $p$ an odd prime. If $p=2$, then
this holds for $n\ne 1,3,7$. 
\end{theorem}

The homologically trivial case was completely settled
by Carlsson \cite{Ca1}, \cite{Ca2}.
More recently the result above was
extended to the case $p=2$
and $n=1$ by Yalcin \cite{Y}. The other cases
remain intriguingly open.
A geometrically interesting situation arises when the action
of $G$ permutes a basis for the homology of the product
of spheres. Here we have a stronger bound, expressed in
the following result, due to Adem and Benson \cite{ABe}.

\begin{theorem}
Let $G$ be an elementary abelian $p$-group of rank $r$
acting freely on a finite dimensional CW complex 
$X\simeq (\mathbb S^n)^t$ in such a way that the basis
$u_1,u_2,\dots , u_t$ of $H^n(X,\mathbb F_p)$ corresponding
to the $t$ spheres is permuted by $G$. Then the number
of orbits of $G$ on $\{ u_1,\dots , u_t\}$ is at least $r$.
\end{theorem}

From this we infer that if $G$ is an elementary abelian
$p$ group acting freely
on $X=(\mathbb S^n)^k$ and permuting
the canonical basis, then $r(G)\le dim~H_n(X,\mathbb F_p)^G$.
Note that if $p=2$ and $n\ne 1,3,7$, then every action
satisfies this hypothesis (by a Hopf invariant one argument,
see \cite{A}). 
We should also mention that Conjecture~\ref{conj-A} follows
from a more general conjecture (due to Carlsson,
see \cite{Ca4}) arising from commutative
algebra, namely:

\begin{conjecture}
If $C_*$ is a free, finite
and connected $\mathbb F_pG$ chain complex, 
where $G=(\mathbb Z/p)^r$, then 
$\sum_{i=0}^{\infty}~dim_{\mathbb F_p}~H_i(C)\ge
2^r.$
\end{conjecture}

This has been verified in some instances; when $p=2$ it
is known to hold for $r_2(G)\le 4$ 
(see \cite{Ca3} and \cite{R}); this
was also settled in the general case of $1$--dimensional chain
complexes in \cite{AS}.

Methods from the cohomology of groups are used to
study these problems. To illustrate some of the techniques
we will prove a simple yet useful result which can help us
approach Conjecture~\ref{conj-A}.\footnote{Although 
this result was known to the author for some time,
it should be attributed to \cite{DGZ}, where this method has
been explained in full detail.}
Recall that $BG$ denotes the classifying
space of a group $G$; given any free $G$--action $X$, there
is a classifying map $X/G\to BG$ associated to the action.

\begin{proposition}
Let $G=(\mathbb Z/2)^r$ act freely on a closed 
n-manifold $M$
such that the classifying map $\pi_G: M/G\to BG$
induces a surjection onto the top non--zero
mod 2 cohomology group $H^n(M/G,\mathbb F_2)$. Then the action
does not extend to a free
action of a larger elementary abelian
$2$--group.
\end{proposition}
\begin{proof}
Assume then that $G$ acts as stated in the theorem, but
that the action extends to a free action by an elementary
abelian $2$--group $K$, with 
$[K:G]=2$. Recall that given an index two
subgroup there is a Gysin sequence associated to it in
equivariant cohomology, which is natural with respect
to equivariant maps. Hence we obtain a diagram of Gysin
sequences, where we assume coefficients are in
$\mathbb F_2$:

\[
\xymatrix{
H^n(BG)\ar[r]^{\Tr} \ar[d]^{\pi_G^*} & 
H^n(BK)\ar[r]^{\cup\chi} \ar[d]^{\pi_K^*}& H^{n+1}(BK)
\ar[d]^{\pi_K^*} \\ 
H^n(M/G)\ar[r]^{\Tr} & H^n(M/K)\ar[r]^{\cup\chi}
 & H^{n+1}(M/K) \\ 
}
\]
Here $\chi\in H^1(BK)$ is the cohomology class
associated to the index two subgroup $G\subset K$.
Choose $0\ne x\in H^n(M/K)$, then as $H^{n+1}(M/K)=0$,
there exists a class $0\ne x'\in H^n(M/G)$ such that
$x=\Tr(x')$. Using our hypotheses,
we can find an $\alpha\in H^n(BG)$
with $\Tr(\pi_G^*(\alpha))=x$. However by the commutativity of the
diagram, we have that this expression is the same as
$\pi_K^*\Tr(\alpha)=0$, as the transfer map is zero for
elementary abelian 2--groups. The result follows.
\end{proof}

For example, one can show that the hypotheses of this theorem
hold for actions of $\mathbb Z/2\times\mathbb Z/2$ on a product
of two spheres (see \cite{DGZ} for details) and more generally
for actions of $(\mathbb Z/2)^r$ on $X=(\mathbb S^n)^r$ 
(see \cite{Ca1}).

From this we derive a conjecture which would imply Conjecture~\ref{conj-A},
namely:

\begin{conjecture}
If $G=(\mathbb Z/2)^r$ acts freely on $X=
\mathbb S^{n_1}\times\dots\times \mathbb S^{n_m}$, then 
$\pi_G^*:H^N(G,\mathbb F_2)\to H^N(X/G,\mathbb F_2)$
is surjective, where $N=n_1 +\dots + n_m$.
\end{conjecture}

\begin{remark}
The case of $(\mathbb Z/2)^r$ actions on a product of real
projective spaces is also of some interest.
For free actions, the main result there (see \cite{AY}) is: if $(\mathbb Z/2)^r$ acts freely on
a finite complex $X\simeq (\mathbb RP^n)^m$, then

$$
r\le \left\{ \begin{array}{ccc}
0 & n\equiv 0,2 &\quad \mbox{mod}~ 4\\
m & n\equiv 1 &\quad\mbox{mod}~ 4\\
2m & n\equiv 3 &\quad\mbox{mod}~ 4\end{array}\right. 
$$

The case of non--free actions is worthy of some attention. For
example, we have that if $(\mathbb Z/2)^k$ acts on
$X=(\mathbb S^m)^t$, then there exists an $x\in X$ with
$\rm{dim} G_x\ge k-t$ (this follows from the 
methods used in \cite{Ca1}). In contrast,
actions on real projective spaces behave very differently; for
example we have (see \cite{AY}):

\begin{proposition}
There exists an action of $G=(\mathbb Z/2)^{1249}$ on
$X=(\mathbb RP^{2^{1298}-1})^{50}$, with isotropy subgroups of
rank less than or equal to 50.
\end{proposition}

\noindent This example arises from representations for certain class 2
nilpotent groups, introduced by Ol'shanskii \cite{O}.
\end{remark}

\section{Constructing Free Actions}

Given the results in the previous section, a natural problem
is that of \emph{constructing} free actions on a product of
of $r(G)$ spheres. In his landmark paper \cite{S}, Swan  
showed that if  $r(G)=1$, then $G$ in fact acts freely
on a finite complex $X\simeq\mathbb S^n$. 
The natural extension of this result leads to a difficult
open problem:

\begin{conjecture}\label{conj-B} 
A finite group $G$ acts freely on a finite complex
$X$ which has the homotopy type of a product of $r(G)$
spheres.
\end{conjecture}
This conjecture is due to Benson and Carlson.
In a purely algebraic context, this was completely settled
in \cite{BC}.
The geometric side of the problem is very difficult, as it
amounts to explicitly constructing free group actions on
products of spheres. 

The most basic examples of spheres with actions of finite groups
arise from representation theory. Given a unitary or orthogonal
$G$--representation $V$, the unit vectors $S(V)$ are homeomorphic
to a sphere. In many cases one can choose representations $V_1,V_2,\dots ,
V_r$ such that $G$ acts freely on $S(V_1)\times\dots\times S(V_r)$. 
However, there are some restrictions, as the following result due to
U.Ray (\cite{Ra}) shows.

\begin{proposition}
If $G$ is a finite group acting freely on a product of spheres arising
from representations, then the only possible non--abelian composition
factors of $G$ are the alternating groups $A_5$ and $A_6$.
\end{proposition}

\begin{example}
Let $P$ denote the $p$--Sylow subgroup of $\rm{GL}_3(\mathbb F_p)$;
it is an extra--special group of order $p^3$. It is not hard to show
that $P$ acts freely on a product $S(V_1)\times S(V_2)\times S(V_3)$;
however more is true (see \cite{Ra}): in fact the \emph{minimal} number of 
representations required for a free action on a product is always three.
\end{example}

\begin{example}
In \cite{O}, it was shown that the alternating group $A_4$
cannot act freely on any finite complex 
$X\simeq \mathbb S^n\times \mathbb S^n$, for any $n$. Hence in
particular $A_4$ cannot act freely on any $\mathbb S^n\times \mathbb S^m$
via a product action, as we could produce a free action on a product
of two equidimensional spheres by taking the appropriate joins.
\end{example}

To construct new examples of free group actions we will
build on what we know about groups with periodic 
cohomology. In fact we will consider $G$--CW complexes
such that the isotropy subgroups all satisfy this 
condition. The following result appears in \cite{ASm}.

\begin{theorem}\label{adem-smith}
Let $X$ denote a simply--connected finite $G$--CW complex
such that all of its isotropy subgroups have periodic 
cohomology. Then there exists a finite complex 
$Y\simeq X\times \mathbb S^n$ with a free $G$--action.
\end{theorem}

This result can be proved directly (see \cite{U})
by constructing a spherical
fibration over successive skeleta in $X$. This uses coherence
provided by a ``universal'' periodicity class $\alpha
\in H^*(G,\mathbb Z)$ for actions with periodic isotropy;
classical spherical space forms and the fact that all
obstructions can be killed using fiber joins.
This is based on the methods first developed in \cite{CP}.
The result is also a consequence of a more general theorem
in \cite{AS} about spaces with periodic cohomomology.
This notion (requiring twisted coefficients)
is induced by cup product with a given 
cohomology class.

From the previous result we see that a clear new strategy
emerges for constructing free group actions: to get a free
$G$--action on $X\times \mathbb S^n$ for some $n>0$, it suffices
to build an action of $G$ on $X$ with periodic 
isotropy.\footnote{These techniques also apply 
to infinite groups, and they can be used to show (see \cite{ASm})
that a discrete group $\Gamma$ acts freely and properly
discontinuously on $\mathbb S^n\times\mathbb R^m$ for some
$m,n>0$ if and only if $\Gamma$ is countable and has
periodic cohomology.}
We now look more carefully at finite group actions on spheres; given
$X\simeq \mathbb S^N$ with a $G$--action, we have an Euler class
$\alpha\in H^{N+1}(G,\mathbb Z)$ associated to it.

\begin{definition}
An Euler class $\alpha\in H^{N+1}(G,\mathbb Z)$ is said to be
effective if $H^*(X\times_GEG,\mathbb Z)$ has Krull dimension
less than the rank $r(G)$.
\end{definition}

\begin{question}\label{quest1}
Given a finite group $G$, does there exist an effective
Euler class in $H^*(G,\mathbb Z)$?
\end{question}

If we assume $X$ to be finite dimensional, then we see that
the Euler class $\alpha$ is effective if and only if 
$X^E=\emptyset$ for all $E\subset G$ elementary abelian
subgroups of maximal rank.  

We now focus on the case of a finite $p$--group $P$. It follows from
known results (see \cite{DH}) that if an effective Euler
class exists, then an effective Euler class must arise from
a representation sphere $S(V)$. Let $C\cong\mathbb Z/p$ be a central
subgroup in $P$, and let $V=\rm{Ind}_C^P(\chi)$
where $\chi$ is a non--trivial character for $C$. Then
$X=S(V)$ has an action of $P$ such that $C$ acts freely.
Hence $X^E=\emptyset$ for all $E$ elementary abelian
subgroups of maximal rank and so we obtain an effective
Euler class $\alpha\in H^{2[P:C]}(P,\mathbb Z)$. More generally,
if the center $Z\subset P$ has rank equal to $z$, then we
can find representations $V_1,V_2,\dots , V_z$ of $P$ such
that $Z$ acts freely on the product $S(V_1)\times\dots
\times S(V_z)$. Note that the isotropy for the action of
$P$ will have rank equal to $r(P)-r(Z)$. This indicates
that for any $p$--group we can construct at least part of the
desired free action on a product of spheres ``for free''.
This is analogous to the fact that we always have a regular
sequence of length $r(Z)$ in $H^*(P,\mathbb F_p)$.

\section{Rank Two Groups}

We now consider groups $G$ with $r(G)=2$. If $P$ is a $p$--group
of rank equal to two, then we know that it will act on a sphere
$X=S(V)$ with periodic isotropy. Hence we obtain (see \cite{ASm}):

\begin{theorem}
A finite $p$--group $P$ acts freely on a finite complex
$X\simeq \mathbb S^n\times \mathbb S^m$ if and only if it
does not contain a subgroup isomorphic to $(\mathbb Z/p)^3$.
\end{theorem}

Of course in particular this establishes Conjecture~\ref{conj-B}
for $r(P)=2$. In addition the conjecture has been settled
affirmatively for almost all rank 2 simple groups. Indeed,
we have

\begin{theorem}
If $G$ is a finite simple group with $r(G)=2$
other than $PSL_3(\mathbb F_p)$, then it acts freely on a 
finite complex $X\simeq\mathbb S^n\times\mathbb S^m$; moreover
any such action 
must be exotic, i.e. it cannot be a product action.
\end{theorem}

The proof that such an action cannot be a product action is
a consequence of the classification of finite simple groups!
Indeed, one checks every non--abelian simple group contains a copy
of $A_4$, and we know that this group cannot act freely
via a product action.

\begin{example}
Let $G=\rm{SL}_3(\mathbb F_2)$, $|G|=2^3\cdot 3\cdot 7$.
From the character tables in the \emph{Atlas of Finite
Groups} \cite{Conway}, we note that
there exists a $3$--dimensional complex representation
$V$ of $G$ such that $V^E=\{ 0\}$ for every
$E=\mathbb Z/2\times\mathbb Z/2$ in $G$. Hence
$X=S(V)$ is a $G$--space with periodic isotropy, and
so, applying our previous results, $G$ acts freely on a 
finite complex 
$Y\simeq \mathbb S^5\times\mathbb S^m$ for some $m$.
However in this case we can proceed more directly as follows.
The subgroup $G\subset U(3)$ must actually lie in $SU(3)$,
as it is a simple group. Then if we consider the $SU(2)$--fibration
$SU(3)\to SU(3)/SU(2)=\mathbb S^5$, we see that it is $G$--equivariant
and that the $G$--action on the total space is free. Hence our
group acts freely on the total space $E$ of an $SU(2)$--bundle
over $\mathbb S^5$, which we can regard  
as the sphere bundle of a 2-dimensional
complex $G$--vector bundle 
$\tau$ over $\mathbb S^5$. This corresponds (non-equivariantly)
to the 
non--zero element
in $\pi_4(SU(2))\cong\mathbb Z/2$. 
Now the sphere bundle $S(\tau\oplus\tau)$
is a trivial bundle, as any $\mathbb S^7$--bundle over
$\mathbb S^5$ splits as a product (see \cite{St}, page 139).
Hence $G=\rm{SL}_3(\mathbb F_2)$ acts freely and smoothly on the 
manifold\footnote{More generally we have shown that if $G$ is any
finite subgroup of $SU(3)$, then it acts freely and smoothly on
$\mathbb S^5\times\mathbb S^7$; this is a special case of a
result in \cite{ADU}.}
$M=\mathbb S^5\times\mathbb S^7$.
\end{example}

Our techniques reduce the problem of
constructing a free action on a product of two
spheres to constructing an action on a 
single sphere with periodic isotropy.
Equivalently we can construct an action on a sphere such that
the associated Euler class is effective.
Although this is a considerable reduction, there are
still groups which cannot be handled with this
approach. In particular we have a result due
to Unlu \cite{U}:

\begin{proposition}
If $G=PSL_3(\mathbb F_p)$ acts on a finite complex
$X\simeq \mathbb S^n$, then there exists a subgroup
$H\cong\mathbb Z/p\times\mathbb Z/p$ in $G$ such that
$X^H\ne\emptyset$.
\end{proposition}

In fact this result has been recently generalized by
Grodal and Smith to show that if $G$ acts on 
\emph{any} homotopy sphere $X$, then the equivariant
cohomology $H^*(EG\times_GX,\mathbb F_p)$ has Krull
Dimension equal to two, hence is not periodic in the
sense of \cite{ASm} and so there is no hope of
constructing a free $G$ action on a product of two spheres
using the methods above.
In particular this gives an example of a group for which
Question~\ref{quest1} has a negative answer.

Motivated by these considerations and in view of some very
special properties of the group cohomology
$H^*(PSL_3(\mathbb F_3),\mathbb F_3)$, it seems valid to raise the
following

\begin{question}
Does $G=PSL_3(\mathbb F_3)$ in fact act freely on a finite complex
$X\simeq \mathbb S^n\times \mathbb S^m$?
\end{question}

Needless to say it would be rather significant if this
question had a negative answer. 

On the other hand this seems to be a rather
special situation, which does not apply for example
to odd order groups. Work in progress (joint with Grodal)
indicates that every odd order group of rank equal to two
will indeed act freely on a finite 
$X\simeq\mathbb S^n\times\mathbb S^m$.
This requires a detailed analysis of
their $p$--Sylow subgroups and fusion properties.

\section{Geometric Actions}

So far we have not discussed the more
geometric side of contructing
group actions. Ideally we would like a group of
rank equal to $r$ to act freely on an actual product
of $r$ spheres---by this we mean a manifold homeomorphic
to it, or even better, diffeomorphic. Although this
has been settled in the rank one case, very little is
currently known for higher rank groups. Given our
previous results in the homotopy category, a good 
question to raise would be

\begin{question}
If $P$ is a $p$--group with $r(P)=2$, does it act
freely on the manifold $M= \mathbb S^n\times \mathbb S^m$
for some $m,n >0$?
\end{question}

This problem is rather difficult as it naturally leads
into both homotopy--theoretic and surgery related
issues. We can however offer a very recent result in this
direction (joint work with Davis and Unlu \cite{ADU}):

\begin{theorem}
Let $G\subset U(n)$ be a faithful unitary representation
of $G$, such that the induced action of $G$ on $U(n)/U(k)$
is free. Then, if $|G|$ is relatively prime to $(n-1)!$,
the group $G$ will act freely on
$M=\mathbb S^{2n-1}\times\dots\times \mathbb S^{2k+1}$.
\end{theorem}

\begin{corollary}
Let $P$ denote a $p$--group with cyclic center and an
abelian maximal subgroup. Then $P$ acts freely on
$M=\mathbb S^{2p-1}\times\dots\times \mathbb S^{2p-r-1}$
where $r<p$ is the rank of $P$.
\end{corollary}

\begin{example}
We consider the extra--special $p$--group $P$ of order
$p^3$ and exponent $p$. Applying the above, it will
act freely on $M=\mathbb S^{2p-1}\times\mathbb S^{2p-3}$.
\end{example}

The proof of the theorem above requires an interesting 
combination of group theory, homotopy theory and surgery
theory. The basic strategy is to propagate the free action
on the Stiefel manifold to a product of spheres.
The corollary follows from the classification of
those finite
$p$--groups which admit a faithful unitary representation
in $U(p)$ such that the induced action on the complex
Stiefel manifold $U(p)/U(p-r)$ is free.

\end{document}